\input amstex

\documentstyle{amsppt}
\loadbold

\def\<{\left<}										%for inner products
\def\>{\right>}

\nologo

%This is Michael's version of Halmos' tombstone.  
\def\qedd{
\hfill
\vrule height4pt width3pt depth2pt
\vskip .5cm
}

\magnification=\magstephalf
%\magnification=1200

\topmatter
\title
The heat flow of the CCR algebra
\endtitle

\author William Arveson
\endauthor

\affil Department of Mathematics\\
University of California\\Berkeley CA 94720, USA
\endaffil

\date 5 April, 2000
\enddate
\thanks On appointment as a Miller Research 
Professor in the Miller Institute for Basic 
Research in Science.  
Support is also acknowledged from 
NSF grant DMS-9802474
\endthanks
%
%\keywords 
%\endkeywords
%
%\subjclass
%\endsubjclass
%
\abstract 
Let $Pf(x)=-if^\prime(x)$ and 
$Qf(x) = xf(x)$ be the canonical 
operators acting on an appropriate common dense 
domain in $L^2(\Bbb R)$.  The derivations 
$D_P(A)=i(PA - AP)$ and $D_Q(A)=i(QA-AQ)$ act on the $*$-algebra 
$\Cal A$ of all integral operators having smooth kernels of compact 
support, for example, and one may consider the 
noncommutative ``Laplacian" $L=D_P^2 + D_Q^2$ as a 
linear mapping of $\Cal A$ into itself.  

$L$ generates a semigroup of normal completely positive 
linear maps on $\Cal B(L^2(\Bbb R))$, and we establish some basic 
properties of this semigroup and its minimal dilation to 
an $E_0$-semigroup.  In particular, we show that its minimal 
dilation is pure, has no normal invariant states, and in section 3 
we discuss 
the significance of those facts for the interaction theory 
introduced in a previous paper.  

There are similar results for the canonical commutation relations with 
$n$ degrees of freedom, $1\leq n <\infty$.  
\endabstract

%\leftheadtext{Blah}
\rightheadtext{The $CCR$ heat flow}
\endtopmatter

\document

\subhead{1.  Discussion, basic results}
\endsubhead

Consider the canonical operators $P,Q$ acting on an 
appropriate common dense domain in $L^2(\Bbb R)$
$$
\align
P&=\frac{1}{i}\cdot\frac{d}{dx},\\
Q&={\text{multiplication by }}x.  
\endalign
$$
These operators can be used to define 
unbounded derivations (say on the dense 
$*$-algebra $\Cal A$ of all integral operators having kernels 
which are smooth and of compact support) by 
$$
D_P(X) = i(PX-XP), \qquad D_Q(X)=i(QX-XQ), \qquad X\in\Cal A.  
$$
Thinking of these derivations as noncommutative counterparts of 
$\partial/\partial x$ and 
$\partial/\partial y$ we define a ``Laplacian" 
$L: \Cal A\to\Cal A$ by
$$
L = D_P^2 + D_Q^2.  \tag{1.1}
$$

Throughout this paper we will use the term CP semigroup to denote 
a semigroup $\phi=\{\phi_t: t\geq 0\}$ of normal completely positive 
linear maps on the algebra $\Cal B(H)$ of all bounded operators on 
a separable Hilbert space $H$, which preserves the unit 
$\phi_t(\bold 1)=\bold 1$, and which is continuous in the natural 
sense (namely $\<\phi_t(A)\xi,\eta\>$ should be continuous in $t$ for 
fixed $\xi,\eta\in H$ and $A\in\Cal B(H)$).  
The purpose of this section is to exhibit concretely a 
CP semigroup whose 
generator can be identified with the operator mapping $L$ of (1.1)
(see Theorem 1.10).

Let $U_t=e^{itQ}$, $V_t=e^{itP}$ be the two unitary 
groups associated with $Q$, $P$,
$$
U_tf(x) = e^{itx}f(x),\qquad V_tf(x)=f(x+t), \qquad f\in L^2(\Bbb R).  
$$
These two groups satisfy the Canonical Commutation Relations 
$V_tU_s=e^{ist}U_sV_t$ for $s,t\in\Bbb R$.  It is 
more convenient to make use of the CCRs in Weyl's form.  
For every $z=(x,y)\in\Bbb R^2$ the Weyl operator 
$$
W_z = e^{\frac{ixy}{2}}U_xV_y \tag{1.2}
$$
is unitary, it is strongly continuous in $z$, and it satisfies 
the Weyl relations
$$
W_{z_1}W_{z_2}=e^{i\omega(z_1,z_2)}W_{z_1+z_2} \tag{1.3}
$$
where $\omega$ is the symplectic form on $\Bbb R^2$ given by 
$$
\omega((x,y),(x^\prime,y^\prime))=\frac{1}{2}(x^\prime y-x y^\prime).  \tag{1.4}
$$
A strongly continuous mapping $z\mapsto W_z\in\Cal B(H)$ into the 
unitary operators on some Hilbert space $H$ which satisfies 
(1.3) is called a Weyl system.  
It is well known that the Weyl system (1.2) is irreducible, and hence 
the space of all finite linear combinations of the $W_z$ is a unital 
strongly dense $*$-subalgebra of $\Cal B(L^2(\Bbb R))$.  The 
Stone-von Neumann theorem implies that every Weyl system is unitarily 
equivalent to a direct sum of copies of the concrete Weyl system (1.2).  

Proceeding heuristically for a moment, let $D_P$ and $D_Q$ be 
the derivations above.  After formally differenting the 
relation $V_tU_s=e^{ist}U_sV_t$ we find that 
$$
\alignat 2
D_P(U_x)&=ixU_x, &\qquad D_P(V_y)&=0,\\
D_Q(U_x)&=0, &\qquad D_Q(V_y)&=iyU_y,
\endalignat
$$
hence the action of $L=D_P^2+D_Q^2$ on the Weyl 
system (1.2) is given by 
$$
L(W_z)=-(x^2+y^2)W_z=-|z|^2W_z, \qquad z=(x,y)\in\Bbb R^2.  
$$
After formally exponentiating we find that for $t\geq 0$ 
the operator mapping $\phi_t=\exp(tL)$ for $t\geq 0$ can 
be expected to satisfy 
$$
\phi_t(W_z)=e^{-t|z|^2}W_z, \qquad z\in\Bbb R^2, t\geq 0.  
\tag{1.5}
$$

\remark{Remarks}
A number of authors have considered completely positive semigroups 
defined on a Weyl system by formulas such as (1.5), using techniques 
similar to those of Proposition 1.7 below (see pp. 128-129 of \cite{3},
or \cite{6} for two notable examples).  We include a full discussion 
of these basic issues since in section 3 we require details of the 
construction that are not easily found in the literature.  

We also remark that virtually all of the results below have straightforward
generalizations to the case in which $P,Q$ are 
replaced with the canonical operators $P_1,\dots,P_n$, $Q_1,\dots,Q_n$
associated with $n$ degrees of freedom.  Indeed, the generalization 
amounts to little more than a reinterpretation of notation.  
On the other hand, while 
Proposition 1.7 remains valid (with the same proof) 
for infinitely many degrees of 
freedom (c.f. \cite{3} {\it loc cit}), we do not know 
if that is the case for the more precise results of section 3.    
\endremark

\vskip0.1truein

In order to define the CCR heat flow rigorously we take (1.5) as our 
starting point and deduce the existence of the 
semigroup and its basic properties from the following 
general result.  
Consider the Banach space $M(\Bbb R^2)$ of all 
complex-valued measures $\mu$ on $\Bbb R^2$
having finite total variation $\|\mu\|$.  $M(\Bbb R^2)$ is a commutative 
Banach algebra with unit relative to the usual convolution of measures 
$$
\mu * \nu(S) = \int_{\Bbb R^2\times\Bbb R^2}\chi_S(z+w)\,d\mu(z)\,d\nu(w).  
$$
It will be convenient to define the Fourier transform of a measure 
$\mu\in M(\Bbb R^2)$ in terms of the symplectic form $\omega$ of (1.4)
$$
\hat \mu(\zeta) =\int_{\Bbb R^2}e^{i\omega(\zeta,z)}\,d\mu(z).  \tag{1.6}
$$
\remark{Remark}
While this definition of the Fourier transform differs from the 
usual one, which involves the Euclidean inner product of $\Bbb R^2$ 
$$
\<(x,y),(x^\prime,y^\prime)\>= xx^\prime + yy^\prime
$$
rather than the symplectic form $\omega$, it is equivalent to it 
in a natural way.  Indeed, since 
$\omega$ is nondegenerate there is a unique invertible skew symmetric 
linear operator $\Omega$ on the two dimensional real vector space 
$\Bbb R^2$ satisfying $\omega(z,z^\prime)=\<\Omega z,z^\prime\>$ 
for all $z,z^\prime\in \Bbb R^2$.  Hence one can pass back and 
forth from the usual Fourier transform of a measure to the one above 
by the invertible linear change-of-variables given by composing the 
transformed measure with either $\Omega$ or $\Omega^{-1}=-4\Omega$.   
\endremark

\proclaim{Proposition 1.7}
Let $\{W_z: z\in\Bbb R^2\}$ be an irreducible Weyl system acting 
on a Hilbert space $H$.  For every complex measure 
$\mu\in M(\Bbb R^2)$ there is a 
unique normal completely bounded linear map 
$\phi_\mu: \Cal B(H)\to\Cal B(H)$ satisfying 
$$
\phi_\mu(W_z)=\hat\mu(z)W_z, \qquad z\in \Bbb R^2.  
$$

One has $\phi_\mu \circ\phi_\nu = \phi_{\mu *\nu}$, and 
$\|\phi_\mu\|_{cb}\leq \|\mu\|$, where $\|\psi\|_{cb}$ denotes
the completely bounded norm of an operator mapping $\psi$.  
When $\mu$ is a positive measure $\phi_\mu$ is a completely 
positive map.  
\endproclaim

\demo{proof}
Fix $\mu\in M(\Bbb R^2)$.  The uniqueness of the mapping 
$\phi_\mu$ is apparent from the irreducibility hypothesis on 
the Weyl system, since the 
set of all linear combinations of the $W_z$, $z\in \Bbb R^2$, is 
a unital $*$-algebra which is weak$^*$-dense in $\Cal B(H)$.  

For existence, we exhibit $\phi_\mu(A)$ for $A\in\Cal B(H)$ 
as a weak integral
$$
\phi_\mu(A)=
\int_{\Bbb R^2}W_{2^{-1/2}\zeta}AW_{2^{-1/2}\zeta}^*\,d\mu(\zeta),
\tag{1.8}
$$
namely the operator defined by the bounded sesquilinear form 
on the right of (1.9)
$$
\<\phi_\mu(A)\xi,\eta\>=
\int_{\Bbb R^2}\<W_{2^{-1/2}\zeta}AW_{2^{-1/2}\zeta}^*\xi,\eta\>
\,d\mu(\zeta), \qquad \xi,\eta\in H.  \tag{1.9}
$$
A straightforward estimate shows that 
$\|\phi_\mu(A)\|\leq \|A\|\cdot\|\mu\|$, and after promoting $\phi_\mu$ 
to $n\times n$ matrices over $\Cal B(H)$ a similar estimate shows 
that $\|\phi_\mu\|_{cb}\leq\|\mu\|$.  $\phi_\mu$ is obviously 
completely positive when $\mu$ is a positive measure.  

Formula (1.9), together with a straightforward application 
of the bounded convergence theorem,  implies that 
when $A_1,A_2,\dots$ is a (necessarily bounded) sequence in 
$\Cal B(H)$ which converges weakly to $A$, one has
$$
\lim_{n\to\infty}\<\phi_\mu(A)\xi,\eta\>=\<\phi_\mu(A)\xi,\eta\>,
\qquad \xi,\eta\in H.  
$$
It follows that $\phi_\mu$ is a normal linear map.  

Finally, from the commuation relation (1.3) we find that 
$$
W_{2^{-1/2}\zeta}W_zW_{2^{-1/2}\zeta}^*=
W_{2^{-1/2}\zeta}W_zW_{-2^{-1/2}\zeta}=e^{-\omega(\zeta,z)}W_z,
$$
and hence (1.8) implies that $\phi_\mu(W_z)=\hat\mu(z)W_z$.\qedd
\enddemo

I want to thank Daniel Markiewicz for a suggestion that simplified  
the proof of Proposition 1.7.

\proclaim{Theorem 1.10}
Let $W=\{W_z: z\in\Bbb R^2\}$ be an irreducible Weyl system.  Then there 
is a unique CP semigroup $\phi=\{\phi_t: t\geq 0\}$ satisfying 
$$
\phi_t(W_z)=e^{-t|z|^2}W_z, \qquad z\in\Bbb R^2.  \tag{1.11}
$$
The only bounded normal linear functional $\rho$ for which 
$\rho\circ\phi_t=\rho$ for all $t\geq 0$ is $\rho=0$.  In particular, 
there is no normal state of $\Cal B(H)$ which is invariant under $\phi$.  
\endproclaim

\demo{proof}  For each $t\geq 0$, $u_t(z)=e^{-t|z|^2}$ is a continuous 
function of positive type, which takes the value $1$ at $z=0$.  Thus it 
is the Fourier transform of a unique probability measure $\mu_t\in M(\Bbb R^2)$.  
We will require an explicit formula for the Gaussian measure $\mu_t$ 
later on; but for purposes of this section we require nothing 
more than its existence and uniqueness.  

Since $u_s(z)u_t(z)=u_{s+t}(z)$ for all $z\in\Bbb R^2$ 
it follows that $\mu_s * \mu_t = \mu_{s+t}$.  Hence Proposition 1.7 implies 
that there is a semigroup $\phi=\{\phi_t: t\geq 0\}$ of normal completely 
positive maps on $\Cal B(H)$ which satisfies (1.11).  It is a 
simple matter to check that 
the required continuity of $\phi_t$ in $t$ follows from the continuity of 
the right side of (1.11) in $t$ for fixed $z$.  

Suppose now that $\rho$ is a normal linear functional 
which is invariant under $\phi$.  Then 
for every $z\in\Bbb R^2$ 
and every $t\geq 0$, the definition of $\phi_t$ implies that 
$$
\rho(W_z)=\rho(\phi_t(W_z))=e^{-t|z|^2}\rho(W_z)
$$ 
and for fixed $z\neq0$, the right side 
tends to $0$ as $t\to\infty$.  Hence $\rho(W_z)=0$
for every $z\neq 0$; by strong continuity on the unit 
ball it follows that $\rho(\bold 1) =\omega(W_0)=0$.  
hence $\rho$ vanishes on the irreducible $*$-algebra 
spanned by $W_z$, $z\in\Bbb R^2$ and by normality it 
follows that $\rho=0$.  \qedd
\enddemo

\remark{Remarks}
We point out that while $\phi$ has no normal invariant
states, it does have a normal invariant weight...namely 
the trace, in that 
$$
{\text{trace}}\,(\phi_t(A))={\text{trace}}\,(A)
$$
for every positive operator $A\in\Cal B(L^2(\Bbb R))$ and 
every $t\geq 0$.  One sees this immediately from (1.8).  
It follows that $\phi_t$ leaves the $C^*$-algebra 
$\Cal K$ of all compact operators invariant, 
$\phi_t(\Cal K)\subseteq\Cal K$.  Since $K$ is the 
$C^*$-algebra associated with the canonical commuation 
relations (more precisely, $\Cal K$ is the enveloping 
$C^*$-algebra of the Banach $*$-algebra of 
all Weyl integral operators associated with the 
canonical commutation relations with a finite 
number of degrees of freedom), this justifies viewing 
the semigroup of restrictions 
$\{\phi_t\restriction_\Cal K: t\geq 0\}$ as the 
heat flow of the canonical commutation relations.  

We also remark that one can deduce the existence of other 
CP semigroups along similar lines.  For example, the proof 
of Theorem 1.10 implies that there is a ``Cauchy" semigroup 
$\psi=\{\psi_t:t\geq 0\}$ which is defined uniquely by 
the requirement 
$$
\psi_t(W_z)=e^{-t(|x|+|y|)}W_z, \qquad z=(x,y)\in\Bbb R^2, 
$$
and which has properties similar to those discussed above
for $\phi=\{\phi_t: t\geq 0\}$.  
\endremark

\subhead{2. Harmonic analysis of the commutation relations}
\endsubhead

A classical theorem of Beurling asserts that singletons 
obey spectral synthesis.  More precisely, if $G$ is a locally 
compact abelian group and $f$ is an integrable function 
on $G$ whose Fourier transform vanishes at a point $p$ 
in the dual of $G$, 
then there is a sequence of functions $f_n\in L^1(G)$ such 
that $\|f - f_n\|\to 0$ as $n\to\infty$ and such that 
the Fourier transform of each $f_n$ vanishes identically on 
some open neighborhood $U_n$ of $p$.  The purpose of this 
section is to present a noncommutative version of that 
result, which will be required in section 3.   

Let $\{W_z: z\in\Bbb R^2\}$ be an irreducible Weyl 
system acting on a Hilbert space $H$ (for example, 
one may take the Weyl system (1.2) acting on 
$L^2(\Bbb R^2)$).  For every trace-class operator $A\in\Cal L^1(H)$ 
we consider the following analogue of the 
Fourier transform $\hat A: \hat\Bbb R^2\to\Bbb C$ 
$$
\hat A(z) = {\text{trace}}(AW_z), \qquad z\in \Bbb R^2.  
$$
This transform $A\in\Cal L^1(H)\mapsto \hat A$ shares many features in
common with the commutative Fourier transform.  For example, 
using the concrete realization (1.2), it is quite easy to establish 
a version of the Riemann-Lebesgue lemma
$$
\lim_{|z|\to\infty}\hat A(z) = 0, 
$$
for every $A\in\Cal L^1(H)$.  What we actually require is 
the following analogue of Beurling's theorem, which lies 
somewhat deeper.  

\proclaim{Theorem 2.1}
Let $A\in\Cal L^1(H)$ and let $\zeta\in\Bbb R^2$ be such 
that ${\text{trace}}(AW_\zeta)=0$.  There is a sequence 
$A_n\in\Cal L^1(H)$ and a sequence of open neighborhoods $U_n$ 
of $\zeta$ such that 
$$
{\text{trace}}(A_nW_z)=0, \qquad z\in U_n, 
$$
and such that ${\text{trace}}|A-A_n|\to 0$ as $n\to\infty$.  
\endproclaim

\demo{proof}
By replacing $A$ with $AW_\zeta$ and making obvious use of 
the canonical commutation relations (1.3), we may immediately 
reduce to the case $\zeta=0$.  We find it more convenient to establish 
the dual assertion of Theorem 2.1.  For that, consider the following linear 
subspaces of $\Cal B(H)$
$$
\Cal S_\epsilon = \overline{{\text{span}}}\{W_z: |z|\leq \epsilon\},
\qquad \epsilon>0,
$$
the closure being taken relative to the weak$^*$ topology on 
$\Cal B(H)$.  Obviously the spaces 
$\Cal S_{\epsilon}$ decrease as $\epsilon$ decreases, and the identity 
operator belongs to $\Cal S_\epsilon$ for every $\epsilon>0$.  The 
pre-annihilator of $\Cal S_\epsilon$ is identified with the space of all 
trace-class operators $A$ satisfying 
$$
\hat A(z)={\text{trace}}(AW_z)=0,\qquad |z|\leq \epsilon.  \tag{2.2}
$$
\proclaim{Lemma 2.3} Let $\{W_z: z\in\Bbb R^2\}$ be an arbitrary 
Weyl system acting on a separable Hilbert space $H$.  Then 
$
\cap\{\Cal S_\epsilon: \epsilon>0\} = \Bbb C\cdot\bold 1.  
$
\endproclaim

\demo{proof of Lemma 2.3}  Let $\Cal S_0$ denote the intersection
$\cap\{\Cal S_\epsilon: \epsilon>0\}$.  
We have already remarked that the inclusion $\supseteq$ is 
obvious.  For the opposite one, consider the von Neumann algebra $\Cal M$ 
generated by $\{W_z: z\in\Bbb R^2\}$.  $M$ is a factor (of type 
$I_\infty$) because of the Stone-von Neumann theorem.  We will show 
that $\Cal S_0$ is contained in the center of $\Cal M$.  

For that, choose $T\in\Cal S_0$ and consider the operator-valued 
function $z\mapsto W_zTW_z^*$.  We have to show that this function 
is constant; equivalently, we will show that for fixed $\xi$ and 
$\eta$ in $H$, the function 
$$
z\in\Bbb R^2\mapsto \<W_zTW_z^*\xi,\eta\> \tag{2.4}
$$
is constant.  Since the function of (2.4) is bounded and continuous, 
it suffices to show that its spectrum (in the sense of spectral 
synthesis for functions in $L^\infty(\Bbb R^2)$) 
is the singleton $\{0\}$: this is the dual formulation of 
Beurling's theorem cited above.  Thus we have to show that 
for every function $f\in L^1(\Bbb R^2)$ whose Fourier transform 
$$
\hat f(\zeta) = \int_{\Bbb R^2}e^{i\omega(z,\zeta)}f(z)\,dz 
$$
vanishes throughout a neighborhood of the origin $\zeta=0$, we have 
$$
\int_{\Bbb R^2}f(z)\<W_zTW_z\xi,\eta\>\,dz = 0.  \tag{2.5}
$$

Fix such an $f\in L^1(\Bbb R^2)$
and choose $\epsilon>0$ 
small enough so that $\hat f(\zeta)=0$ for all $\zeta$ 
satisfying $|\zeta|\leq \epsilon$.  Since the linear functional 
$$
X\in\Cal B(H)\mapsto \int_{\Bbb R^2}f(z)\<W_zXW_z^*\xi,\eta\>\,dz
$$
is weak$^*$-continuous and $T$ belongs to the 
weak$^*$-closed linear span of operators of the form 
$W_\zeta$ with $|\zeta|\leq \epsilon$, to prove (2.5) it suffices 
to show that for every $\zeta$ with $|\zeta|\leq \epsilon$ we have 
$$
\int_{\Bbb R^2}f(z)\<W_zW_\zeta W_z^*\xi,\eta\>\, dz = 0.  \tag{2.6}
$$
Using the canonical commutation relations we can write 
$$
W_zW_\zeta W_z^* = e^{i\omega(z,\zeta)}W_{z+\zeta}W_{-z}=
e^{\omega(\zeta,-z)}W_\zeta=e^{i\omega(z,\zeta)}W_\zeta.  
$$
Hence the left side of (2.6) becomes 
$$
\int_{\Bbb R^2}f(z)e^{i\omega(z,\zeta)}\<W_\zeta\xi,\eta\>\, dz=
\hat f(\zeta)\<W_\zeta\xi,\eta\>,   
$$
and the latter term vanishes because $\hat f(\zeta)=0$ 
for $|\zeta|\leq \epsilon$.  \qedd
\enddemo

To complete the proof of Theorem 2.1, choose an operator 
$A\in\Cal L^1(H)$ satisfying 
$$
\hat A(0)={\text{trace}}(A)=0, 
$$
and consider the linear functional $\rho$ defined on 
$\Cal B(H)$ by $\rho(T)={\text{trace}}(AT)$.  $\rho$ 
obviously vanishes on $\Bbb C\cdot\bold 1$.  The linear 
spaces $\Cal S_\epsilon$ are weak$^*$-closed and they decrease to 
$\Cal S_0=\Bbb C\cdot \bold 1$ as $\epsilon$ decreases to 
$0$, by Lemma 2.3.  Since $\rho$ is weak$^*$-continuous we 
must have 
$$
\lim_{\epsilon\to 0}\|\rho\restriction_{\Cal S_\epsilon}\|
= \|\rho\restriction_{\Bbb C\cdot \bold 1}\| = 0.  
$$
Thus we can choose a sequence $\epsilon_n\downarrow 0$ so 
that $\|\rho\restriction_{\Cal S_{\epsilon_n}}\|\leq 1/n$ 
for every $n=1,2,\dots$.  We have already pointed out that 
the pre-annihilator of $\Cal S_{\epsilon_n}$ is identified 
with all trace class operators $B$ satisfying 
$$
\hat B(z) = {\text{trace}}(BW_z)=0,\qquad |z|\leq \epsilon_n.  \tag{2.7}
$$
Since 
$\|\rho\restriction_{\Cal S_{\epsilon_n}}\|$ is the trace-norm 
distance from $A$ to the pre-annihilator of $\Cal S_{\epsilon_n}$, 
we conclude that there is a sequence of operators $B_n\in\Cal L^1(H)$
which satisfy ${\text{trace}}(B_nW_z)=0$ for $|z|\leq \epsilon_n$, 
such that ${\text{trace}}|A-B_n|\leq 2/n$, as asserted.  \qedd
\enddemo

\subhead{3.  Purity and dilation theory}
\endsubhead

An $E_0$-semigroup is a CP semigroup 
$\alpha=\{\alpha_t: t\geq 0\}$, acting on $\Cal B(H)$, 
such the the individual maps 
are endomorphisms, $\alpha_t(AB) = \alpha_t(A)\alpha_t(B)$, 
$A,B\in\Cal B(H)$.  An $E_0$-semigroup $\alpha$ is called {\it pure}
if its ``tail" von Neumann algebra is trivial, 
$$
\bigcap_{t\geq 0}\alpha_t(\Cal B(H))=\Bbb C\cdot\bold 1.  \tag{3.1}
$$
It is known that an $E_0$-semigroup is pure iff for any pair of normal 
states $\rho_1,\rho_2$ of $\Cal B(H)$ we have 
$$
\lim_{t\to\infty}\|\rho_1\circ\alpha_t - \rho_2\circ\alpha_t\|=0  \tag{3.2}
$$
see \cite{1}.  

If a pure $E_0$-semigroup $\alpha$ has a normal invariant state $\omega$, 
then the characterization (3.2) implies that $\omega$ must be an {\it absorbing} 
state in the sense that for every normal state 
$\rho$ of $\Cal B(H)$ one has 
$$
\lim_{t\to\infty}\|\rho\circ\alpha_t - \omega\| = 0.  \tag{3.3}
$$
Conversely, if for an arbitrary $E_0$-semigroup $\alpha$ there is 
a state $\omega$ of $\Cal B(H)$ which is absorbing in the sense that 
(3.3) is satisfied for every normal state $\rho$ of $\Cal B(H)$, then 
$\omega$ must be a normal invariant state, and thus by 
(3.2) $\alpha$ must be a pure $E_0$-semigroup.  

In the theory of interactions worked out in \cite{2}, pure $E_0$-semigroups 
occupy a central position, especially those for which there is a normal 
invariant (and therefore absorbing) state.  A natural question 
that emerges from the theory of interactions is 
whether or not every pure $E_0$-semigroup must have a normal invariant 
state.  Now since the state space of $\Cal B(H)$ is weak$^*$-compact, a routine 
application of the Markov-Kakutani fixed point theorem shows shows that 
every $E_0$-semigroup must have invariant states; but invariant states 
obtained by such methods need not be normal.  In this section we exhibit 
a concrete $E_0$-semigroup which is pure but which has no {\it normal} 
invariant states.  This is a result which was asserted (without proof) in 
\cite{2}.  This $E_0$-semigroup is obtained from the CP semigroup of 
Theorem 1.10 by a dilation procedure.  

In order that the minimal dilation of a CP semigroup to an 
$E_0$-semigroup should satisfy (3.1), it is necessary and sufficient that the 
CP semigroup should satisfy property (3.2) (see Proposition 3.5).  
Thus we generalize the definition of pure $E_0$-semigroup as follows.  

\proclaim{Definition 3.4}
A CP semigroup
$\phi$ acting on $\Cal B(H)$ is called pure if for every pair of normal 
states $\rho_1$, $\rho_2$ of $\Cal B(H)$ we have 
$$
\lim_{t\to\infty}\|\rho_1\circ\phi_t-\rho_2\circ\phi_t\|=0.  
$$
\endproclaim

\proclaim{Proposition 3.5}  
Let $\phi=\{\phi_t: t\geq 0\}$ be a pure CP semigroup 
which has no normal invariant state, 
and let $\alpha$ be its minimal dilation 
to an $E_0$-semigroup.  Then $\alpha$ 
satisfies (3.1) and has no normal invariant state.  
\endproclaim
\demo{proof}The proof is straightforward, 
but we require results from \cite{1}.  
Proposition 2.4 of \cite{1} implies that 
$\alpha$ satisfies (3.1).  

To see that $\alpha$ has no normal invariant state, 
we can assume that 
$\alpha$ acts on $\Cal B(H)$ for some Hilbert space 
$H$ and that there is a 
closed subspace $K\subseteq H$ such that $\phi$ 
is the compression of $\alpha$ onto 
$\Cal B(K)=P\Cal B(H)P$, $P$ denoting the projection
of $H$ onto $K$.  
We have $\alpha_t(P)\uparrow \bold 1$ because $\alpha$ 
is minimal over $P$.  So if $\omega$ is any normal 
state of $\Cal B(H)$ which is invariant under 
$\alpha$ then we have 
$$
\omega(P) =\lim_{t\to\infty}\omega(\alpha_t(P))=\omega(\bold 1)= 1.  
$$
Thus the restriction of $\omega$ to $\Cal B(K)=P\Cal B(H)P$ 
defines a normal $\phi$-invariant state on $\Cal B(K)$, contradicting
the hypothesis on $\phi$. \qedd
\enddemo

In the remainder of this section we show that the CP semigroup
defined in Theorem 1.10 is pure.  Once that is established, 
Proposition 3.5 implies that its minimal 
dilation is an $E_0$-semigroup with properties asserted in the 
discussion above.  

\proclaim{Theorem 3.6}
The CP semigroup $\phi$ defined in (1.11) is pure.  
\endproclaim

Before giving the proof, we require 

\proclaim{Lemma 3.7}
For each $t>0$ let $\mu_t$ be the Gaussian measure 
on $\Bbb R^2$ whose Fourier transform (1.6) is given by
$$
\hat\mu_t(z)=e^{-t|z|^2}, \qquad z\in\Bbb R^2, 
$$
and choose $\delta>0$.  
There is a family $\nu_t$, $t>0$, of probability measures on 
$\Bbb R^2$ such that 
\roster
\item"{(i)}" $\hat\nu_t(z)=0$ for all $|z|\geq \delta$ and every t, and 
\item"{(ii)}" $\lim_{t\to\infty}\|\mu_t-\nu_t\|=0$, $\|\cdot\|$ denoting 
the norm of the measure algebra $M(\Bbb R^2)$.  
\endroster
\endproclaim

\demo{proof}
For $t>0$, $\mu_t$ is given by $d\mu_t=u_t(x,y)\,dx\,dy$, where 
$u_t$ is the density
$$
u_t(x,y)=\frac{1}{\pi t}e^{-\frac{x^2+y^2}{4t}}.  
$$
Let $f_t=\sqrt{u_t}$.  $f_t$ belongs to $L^2(\Bbb R^2)$ and 
$\|f_t\|_2=1$.  Choose a function $g\in L^1(\Bbb R^2)$ whose 
Fourier transform 
$$
\hat g(\zeta)=\int_{\Bbb R^2}e^{i\omega(\zeta,z)}g(z)\,dz
$$
satisfies $0\leq \hat g(\zeta)\leq 1$ for all $\zeta$, and 
$$
\hat g(\zeta) = 
\cases
1, &{\text{for }} 0\leq |\zeta|\leq \delta/4\\ 
0, &{\text{for }} |\zeta|\geq \delta/2.  
\endcases
$$
Consider the convolution $g*f_t\in L^2(\Bbb R^2)$ and the 
positive measure 
$$
d\nu_t = |g*f_t|^2\,dx\,dy.  
$$
$\nu_t$ is obviously a positive finite measure.  We claim first that 
the Fourier transform of $\nu_t$ lives in the disk $|\zeta|\leq\delta$.  
Indeed, letting $U_\zeta$ (resp. $T_\zeta$) be the unitary operator 
on $L^2(\Bbb R^2)$ (resp. $L^2(\hat\Bbb R^2)$) given by 
$$
U_\zeta F(z)=e^{i\omega(\zeta,z)}F(z), \qquad 
T_\zeta G(w) = G(w+\zeta), 
$$
we have by the Plancherel theorem 
$$
\align
\hat\nu_t(\zeta)&=\<U_\zeta(g*f_t),g*f_t\>_{L^2(\Bbb R^2)} =
\<T_\zeta(\hat g\hat f_t),\hat g\hat f_t\>_{L^2(\hat\Bbb R^2)} \\&= 
\int_{\hat\Bbb R^2}(\hat g\hat f_t)(w+\zeta)
\overline{(\hat g\hat f_t)}(w)\,dw.  
\endalign
$$
When $|\zeta|\geq\delta$ the integrand on the right vanishes 
identically in $w$ because $\hat g\hat f_t$ is supported in the disk 
of radius $\delta/2$.  Hence $\hat\nu_t(\zeta)=0$ for $|\zeta|\geq \delta$.  

To establish property (ii), it is enough to show that 
$$
\lim_{t\to\infty}\|f_t-g*f_t\|_2=0, \tag{3.8}
$$
since by the Schwarz inequality 
$$
\align
\|\mu_t-\nu_t\|&=\int_{\Bbb R^2}|f_t^2-|g*f_t|^2|\,dz \leq
\int_{\Bbb R^2}|f_t-g*f_t|\cdot|f_t+|g*f_t||\,dz
\\&\leq\|f_t-g*f_t\|_2\cdot\|f_t+|g*f_t|\|_2
\leq \|f_t-g*f_t\|_2(\|f_t\|_2+\|g*f_t\|_2).  
\endalign
$$

To establish (3.8), we use the Plancherel theorem again to write 
$$
\int_{\Bbb R^2}|f_t(z)-g*f_t(z)|^2\,dz=
\int_{\hat\Bbb R^2}|\hat f(\zeta)-\hat g(\zeta)\hat f_t(\zeta)|^2\,d\zeta
=\int_{\Bbb R^2}|1-\hat g(\zeta)|^2\cdot|\hat f_t|^2\,d\zeta.  
$$
The function $|1-\hat g(\zeta)|$ is bounded above by $1$ and it 
vanishes throughout the disk $0\leq|\zeta|\leq\delta/4$.  Hence the term 
on the right is dominated by 
$$
\int_{\{|\zeta|\geq\delta/4\}}|\hat f_t(\zeta)|^2\,d\zeta.  \tag{3.9}
$$
In order to estimate the integral (3.9) we require the explicit formula 
$$
f_t(x,y)=\sqrt{u_t(x,y)}=\frac{1}{\sqrt{\pi t}}e^{-\frac{x^2+y^2}{8t}}.  
$$
The Fourier transform of $f_t$ has the form 
$$
\hat f_t(\zeta)=K\sqrt{t}e^{-2t|\zeta|^2}
$$
where $K$ is a positive constant,  hence (3.9) evaluates to 
$$
K^2 t \int_{\{|\zeta|\geq\delta/4\}}e^{-4t|\zeta|^2}\,d\zeta 
=K^2\int_{S_t}e^{-4(u^2+v^2)}\,du\,dv,   
$$
where $S_t=\{(u,v): \sqrt{u^2+v^2}\geq(\delta/4)\sqrt{t}\}$.  As 
$t\to\infty$ the sets $S_t$ decrease to $\emptyset$, hence the right 
side of the previous expression tends to $0$, and (3.8) is proved.  

The positive measures $\nu_t$ are not necessarily probability measures, 
but in view of the established property (ii), 
$\nu_t(\Bbb R^2)$ must be arbitrarily 
close to $\mu_t(\Bbb R^2)=1$ when $t$ is large.  Hence we  can rescale 
$\nu_t$ in an obvious way to achive $\nu_t(\Bbb R^2)=1$ for all $t>0$ 
as well as the properties (i) and (ii) of Lemma 3.7.  
\qedd
\enddemo

\demo{proof of Theorem 3.6}
Let $W_z$, $z\in\Bbb R^2$ be an irreducible Weyl system 
acting on a Hilbert space $H$, and let 
$\phi=\{\phi_t: t\geq 0\}$ be the CP semigroup defined by 
the condition
$$
\phi_t(W_z)=e^{-t|z|^2}W_z,\qquad z\in\Bbb R^2.  
$$
Choose a pair of normal states $\rho_1$, $\rho_2$ on $\Cal B(H)$, 
and consider their difference $\omega=\rho_1-\rho_2$.  We have to 
show that 
$$
\lim_{t\to\infty}\|\omega\circ\phi_t\| = 0.  \tag{3.10}
$$
For that, let $A$ be the self-adjoint trace-class operator defined 
by ${\text{trace}}(AT)=\omega(T)$, $T\in\Cal B(H)$ and choose 
$\epsilon>0$.  $A$ has trace 
zero, so by Theorem 2.1, we can find a 
self-adjoint trace-class 
operator $A_0$ such that ${\text{trace}}(A_0W_z)=0$ for every $z$ in some 
neighborhood $U$ of $z=0$, and ${\text{trace}}|A-A_0|\leq \epsilon$.  It follows 
that the normal linear functional
$\omega_0(T)={\text{trace}}(A_0T)$ satisfies $\|\omega-\omega_0\|\leq \epsilon$
and $\omega_0(W_z)=0$ for $z\in U$.  

By Lemma 3.5 we can find probability measures $\nu_t$, $t>0$, such 
that $\hat\nu_t(z)$ vanishes for $z\notin U$ and $\|\mu_t-\nu_t\|$ 
tends to $0$ as $t\to\infty$.  For each $t>0$
let $\psi_t$ be the 
completely positive map defined by Proposition 1.7,
$$
\psi_t(W_z)=\hat\nu_t(z)W_z, \qquad z\in\Bbb R^2.  
$$
In order to prove (3.10) we decompose the linear 
functional $\omega\circ\phi_t$ into a sum of three terms as follows
$$
\omega\circ\phi_t = (\omega-\omega_0)\circ\phi_t + 
\omega_0\circ(\phi_t-\psi_t) + \omega_0\circ\psi_t.  \tag{3.11}
$$
The third term on the right of (3.11) is zero because for every $z\in\Bbb R^2$ 
we have 
$$
\omega_0(\psi_t(W_z))=\omega_0(\hat\nu_t(z)W_z)=\hat\nu_t(z)\omega_0(W_z)=0,
$$
since $\hat\nu_t(z)$ vanishes when $z\notin U$ and $\omega_0(W_z)$ 
vanishes when $z\in U$ (recall that the linear span of the $W_z$ for 
$z\in\Bbb R^2$ is a strongly dense $*$-subalgebra of $\Cal B(H)$).  The 
first term on the right of (3.11) is estimated for arbitrary $t$ by 
$$
\|(\omega-\omega_0)\circ\phi_t\|\leq \|\omega-\omega_0\|\leq \epsilon.  
$$
In order to estimate the second term, note that 
$$
\|\phi_t-\psi_t\|\leq \|\mu_t-\nu_t\| \tag{3.12}
$$
for every $t>0$.  Indeed, considering the measure 
$\sigma_t=\mu_t-\nu_t\in M(\Bbb R^2)$, we can write 
$$
\phi_t(W_z)-\psi_t(W_z)=\hat\mu_t(z)W_z-\hat\nu_t(z)W_z=\hat\sigma_t(z)W_z.  
$$
It follows from Proposition 1.7 that the completely bounded norm of 
the operator mapping $\phi_t-\psi_t$ is at most $\|\sigma_t\|=\|\mu_t-\nu_t\|$, 
hence (3.12).  

From (3.11) and these estimates we may conclude that 
$$
\limsup_{t\to\infty}\|\omega\circ\phi_t\|\leq 
\epsilon + \lim_{t\to\infty}\|\mu_t-\nu_t\| + 0=\epsilon.  
$$
Since $\epsilon$ is arbitrary the limit (3.10) is proved.  
\qedd
\enddemo

\end
\Refs
%
\ref\no 1\by Arveson, W. \paper Pure $E_0$-semigroups and 
absorbing states
\jour Comm. Math. Phys. \vol 187 \yr 1997
\pages 19--43
\endref

\ref\no 2\bysame\paper Interactions in noncommutative dynamics
\jour Comm. Math. Phys. \paperinfo to appear
\endref

\ref\no 3\by Davies, E. B.\book Quantum theory of open systems
\publ Academic Press\publaddr London\yr 1976 
\endref

\ref \no 4\bysame\paper Generators of dynamical semigroups
\jour J. Funct. Anal. \vol 34\pages 421--432\yr 1979
\endref

\ref\no 5\by Evans, D.\paper Quantum dynamical semigroups, symmetry
groups, and locality \jour Acta Appl. Math.\vol 2\yr 1984
\pages 333--352
\endref

\ref\no 6 \by Evans, D. and Lewis, J. T. \paper Some semigroups 
of completely positive maps on the $CCR$ algebra
\jour J. Funct. Anal.\vol 26 \yr 1977 \pages 369--377
\endref


\ref\no 7 \bysame \paper Dilations
of irreversible evolutions in algebraic quantum theory
\jour Comm. Dubl. Inst. Adv. Studies, Ser A\vol 24\yr 1977
\endref

\ref\no 8\by Gorini, V., Kossakowski, A. and Sudarshan, E. C. G.
\paper Completely positive semigroups on $N$-level systems
\jour J. Math. Phys.\vol 17\yr 1976\pages 821--825
\endref

\ref\no 9\by Hudson, R. L. and Parthasarathy, K. R. \paper
Stochastic dilations of uniformly continuous completely positive
semigroups \jour Acta Appl. Math. \vol 2\pages 353--378\yr 1984
\endref

\ref\no 10\by Lindblad, G.\paper On the generators of quantum 
dynamical semigroups
\jour Comm. Math. Phys.\vol 48\yr 1976\pages 119
\endref

\ref\no 11\by Mohari, A., Sinha, Kalyan B. \paper Stochastic 
dilation of minimal quantum dynamical semigroups \jour
Proc. Ind. Acad. Sci. \vol 102\yr 1992\pages 159--173
\endref



\endRefs
